\documentclass[20pt]{article}

\input epsf
\usepackage{amssymb}
\usepackage{amsfonts}
\usepackage{amsmath}
\usepackage{mathrsfs}
\usepackage{amsthm}

\numberwithin{equation}{section}
\newtheorem{theorem}{Theorem}[section]
\newtheorem{lemma}{Lemma}[section]

\newtheorem{remark}{Remark}[section]
\newtheorem{example}{Example}[section]

 \newcommand{\R}{\mathbb{R}}
\newcommand{\C}{\mathbb{C}}
\newcommand{\Z}{\mathbb{Z}}

\begin{document}

 \baselineskip 20pt

\textwidth=146truemm
\textheight=207truemm

\begin{center}

 {\bf \Large Multivariate Splines and Polytopes \footnote{
The work of this author was supported in part by NSFC grant (10871196) and
National Basic Research Program of China (973 Program  2010CB832702), and was
performed in part at Technische Universit\"at Berlin. }}
\end{center}

\begin{center}
 {Zhiqiang Xu\footnote{Email: xuzq@lsec.cc.ac.cn.}}

{\small LSEC, Institute of Computational Mathematics, Academy of
Mathematics and System Sciences, Chinese
 Academy of Sciences, Beijing, 100091, China }
\end{center}

\begin{abstract}

In this paper, we use multivariate splines to investigate the volume
of polytopes.   We first present an explicit formula for the
multivariate truncated power, which can be considered as a dual
version of the famous Brion's formula for the volume of polytopes.
 We also prove that the integration of polynomials over polytopes
 can  be dealt with by the multivariate truncated power. Moreover,
 we show that the volume of the cube slicing can be considered as the
maximum value of the box spline. Based on this connection, we give a
simple proof for Good's conjecture, which has been settled before by
 probability methods.

\end{abstract}

 \vspace{0.5cm}

 \baselineskip 20pt
\section{Introduction}

Box splines and  multivariate truncated powers were first introduced in
\cite{deboor2} and \cite{dahemenpower}, respectively. They have a wide range of
important and varied applications in numerical analysis and approximation
theory. From the point of  view of discrete geometry,  box splines and
multivariate truncated powers are closely related to the volume of cube slicing
and the volume of polytopes, respectively. However,  people working in the
discrete geometry do not seem to be fully aware of the means of multivariate
splines. The aim of this paper is to recast some problems related to the
computation of volumes of polytopes and solve them by multivariate splines. We
believe that the means of multivariate splines  sheds some light on problems
concerning polytopes.

The main results in this paper are  as follows. Some of the main results
concern the volumes of polytopes.  The exact computation of the volume of a
polytope $P$ is an important and difficult problem which has close ties to
various mathematical areas. Brion's formula in continuous form (see
\cite{beckbook,brion}), which is well known in discrete geometry, gives an
explicit formula for the volume of polytopes. However, the formula lyings on
the vertex representation of polytopes and  requires the generators for each
vertex cone. Based on the multivariate exponential truncated power \cite{ron},
we give an explicit formula for the multivariate truncated power, which can be
regarded as a dual version of Brion's formula. As another effort,  Lasserre
\cite{lasserre} gave a recursive formula for  the volume of polytopes, which
has become a popular method  today. We  re-prove the formula by an iterative
formula for the multivariate truncated power \cite{numeircalT}, which was
presented by Micchelli.

Another problem addressed is the integration of  continuous functions over
polytopes. This problem has important applications.  For example, in most
finite element integration methods, the domain of integration is decomposed
into polytopes. Hence the integration of real functions over  polytopes is
always required.  In \cite{integ} an exact formula for the integration of
polynomials over simplices is presented. An iterative formula for computing the
integration over polytopes is given in \cite{integ1}, and the formula is
extended in \cite{xu} by F-splines. We shall show that integration of
polynomials over polytopes can be dealt with via the multivariate truncated
power leading  to an explicit formula for the integral of polynomials over
polytopes. As continuous functions on a compact set can be uniformly
approximated by polynomials, this result provides an approximate formula for
integrating continuous functions over polytopes. Moreover, this result also
shows that one can compute the integrals of polynomials over polytopes by
calculating the volumes of polytopes.

The volume of  cube slicing is another active research topic in discrete
geometry (see \cite{zong}). Suppose that $Q_n:=[0,1)^n$ is the unit cube in
$\R^n$ and ${H}$ is an $(n-1)$-dimensional hyperplane of $\R^n$ through of its
center. According to \cite{hensley}, Good conjectured that  ${\rm vol}({H}\cap
Q_n)\geq 1$.  Hensley (unexpectedly) introduced a probabilistic method for the
study of this conjecture and finally solved it \cite{hensley}. In fact, we show
that the conjecture can be reformulated as the following box spline problem:
\begin{equation}\label{eq:maxbox1}
\max_x B(x| (a_1,\ldots,a_n)) \geq \frac{1}{\sqrt{\, \sum_{i=1}^n
a_i^2}},
\end{equation}
where $B(x|(a_1,\cdots,a_n))$ is a univariate box spline (cf. Section 3), and
$a_i$ are  positive real numbers for $1\leq i\leq n$. Based on the explicit
formula for the Fourier transform of the box spline, we give a simple proof of
(\ref{eq:maxbox1}). Hence we present a spline method for proving the
conjecture.

The problem of computing the volume of the intersection of $Q_n$ and an
$(n-1)$-hyperplane is also interesting and it can be traced back to P\'{o}lya's
thesis. In \cite{slice}, the authors derived a formula for the volumes of such
domains using combinatorial methods. Using box splines, we give an explicit
formula for the volume of convex bodies obtained by intersecting $Q_n$ and a
$j$-hyperplane, where $j$ is a positive integer $<\, n$. The formula in
\cite{slice} may be considered as a special case of ours.

The paper is organized as follows. After recalling some definitions
and notations in Section 2, we show (Section 3) the connection
between the multivariate truncated power and the volume of
polytopes. In Section 4, we transform  the integration  of
polynomials over polytopes to a problem concerning
 the multivariate truncated power. In Section 5, we
investigate the volume of cube slicing using  box splines. Finally,
Section 6 illustrates the application of the formulas given in this
paper with some examples.

\section{Definitions and Notations }

 A convex polytope $P$ is the convex hull of a finite set of points in
$\R^d$. Throughout this paper, we shall omit the qualifier ``convex" since we
confine our discussion to such polytopes. We also assume that the affine hull
of $P$ contains the origin.  Moreover, we use $d$-polytope to mean a
$d$-dimensional polytope. When the polytope $P$ is defined as the convex hull
of a finite set
 of points in $\R^d$, the finite set is called as a vertex
representation  or simply ${\mathcal V}$-representation of $P$,
while, if $P$ is defined as $\{x\in \R_+^n|\,\, Mx=b\}$ for some
$s\times n$ matrix $M$ and $s$-vector $b$, then the pair $(M,b)$ is
called a half space representation or simply ${\mathcal
H}$-representation. For a vertex $v$ of $P$, we define the vertex
cone of $v$ as the smallest cone with vertex $v$ that contains $P$.
If $P\subset \R^n$ is a $d$-polytope, then let ${\rm vol}_n(P)$
denote the $d$-dimensional volume of $P$ in $\R^n$.  For a rational
polytope $P$, i.e., a polytope whose vertices have rational
coordinates, let $\R_P$ denote the space that is spanned by the
vertex vectors of $P$. The lattice points in $\R_P$ form an Abelian
group of rank $d$, i.e., $\R_P\cap \Z^d$ is isomorphic to ${\Z}^d$.
Hence there exists an invertible affine linear transformation
$T:\R_P \rightarrow \R^d$ satisfying $T(\R_P\cap \Z^n)=\Z^d$. The
{\it relative volume} of $P$, denoted as ${\rm vol}(P)$, is just the
$d$-dimensional volume of the image $T(P)\subset {\R}^d$. For more
detailed information about the relative volume, the reader is
referred to \cite{stanley}.

Throughout this paper, $\Z_+$ and $\R_+$ denote the non-negative
integer and non-negative real sets, respectively. Given a set $D$,
let $\chi_{D}(x) =1$ if $ x\in D$, otherwise let $\chi_{D}(x) =0$.
Elements of $\R^s$ can be regarded as row or column $s$-vectors
according to circumstances. Let $M$ be an $s\times n$ matrix. Then
$M$ can be considered as a multiset of its columns. The cone spanned
by $M$, denoted by ${\rm cone}(M)$, is the set $ \{ \sum_{m\in M}a_m
m|\,\, a_m\geq 0 \mbox{ for all\ } m \}$. Moreover, we set
$[[M)):=\{\sum_{m\in M}a_mm|\,\, 0\leq a_m< 1, \mbox{ for all } m\in
M\}$. Furthermore, we use $\#A$ to denote the cardinality of the
finite set $A$. `` $\cdot$ " stands for scalar product and $\|\cdot
\|$ for the Euclidean norm. $E_{s\times s}$ denotes the $s\times s$
identity matrix. As a final piece of notation, $A^{-T}:=(A^{-1})^T$
when $A$ is an invertible matrix.

\section{Multivariate truncated powers and the volume of polytopes}

Let $M$ be an $s\times n$ real matrix with ${\rm rank}(M)=s$. Recall that $M$
is also viewed as the multiset of its column vectors. Throughout this section
we always assume that the interior of the convex hull of $M$ does not contain
the origin. The {\it multivariate truncated power} $T(\cdot |M)$ associated
with $M$, first introduced by  Dahmen \cite{dahemenpower}, is the distribution
given by the rule
\begin{equation}\label{powerdefinition}
\int_{\R^s}\!\! T(x |M)\phi (x) dx = \int _{\R_+^n}\!\! \phi
(M{u})du,\,\,\, \phi \in {\mathscr D}(\R^s),
\end{equation}
where ${\mathscr D}(\R^s)$ is the space of test functions on
${\R}^s$. If we define  $P:=\{{ y}\in \R_+^n|\,\, M{y}=x\}$, then
(see \cite{deboorbook})
\begin{equation}\label{Eq:truncatedvolume}
T(x|M)=\frac{{\rm vol}_{n}(P)}{\sqrt{\, \det(MM^T)}}.
\end{equation}
 Note that
$$
\frac{{\rm vol}_n(P)}{{\rm vol}(P)}=\frac{\sqrt{\,
\det(MM^T)}}{\#\{[[M^T))\cap {\Z}^s\}}
$$
provided $M$ is an integer matrix (see \cite{barv}). Hence we have
\begin{equation}\label{Eq:truncated1}
T(x| M)=\frac{{\rm vol}(P)}{\#\{[[M^T))\cap \Z^s\}}
\end{equation}
 provided $M$ is an integer matrix.
In particular, if $E_{s\times s}\subset M$ then $T(x |M)={\rm
vol}(P)$, since $\#\{[[M^T))\cap \Z^s\}=1$.
 So, the relative volume of polytopes $P$ can be obtained by computing
$T(x|M)$.

  In the following, we shall give an explicit formula for $T(x|M)$.
   We first introduce the {\it multivariate exponential truncated power}
$E_{c}(x|M)$ associated with a complex vector
${c}=(c_1,\ldots,c_n)\in \C^n$ and a matrix $M$. The function $E_{c}(x|M)$ is the
distribution given by the rule (see \cite{ron}):
\begin{equation}\label{eq:etrun}
\int_{\R^s} \!\! E_{c}(x|M)\phi(x) dx=\int_{\R_+^n}\!\!
\exp{(-{c}\cdot u)}\phi(Mu)du,\ \ \ \phi(x)\in {\mathscr D}(\R^s).
\end{equation}
It is convenient for us to index the constants $c_1,\ldots,c_n$ by
an element $m_i\in M$, that is, we set $c_{m_i}:=c_i$ for
$i=1,\ldots,n$. For the submatrix $M'=(m_{i_1},\ldots,m_{i_k})$ in
$M$ we set $c_{M'}:=(c_{{i_1}},\ldots,c_{{i_k}})$ and
$M'/c_{M'}:=(m_{i_1}/c_{i_1},\ldots,m_{i_k}/c_{i_k})$.

We  recall an explicit formula for $E_{c}(\cdot |M)$. In this
formula, we denote, given a square invertible $Y\subset M$,
$\theta_Y:=Y^{-T}c_Y$ and  $\alpha_Y:=\prod_{y\in M\setminus Y}
\left(\theta_Y\cdot y-c_y\right)^{-1}$.

\begin{lemma}\label{le:E}{\rm (\cite{ron})}
\begin{equation}\label{eq:lee}
E_{c}(x|M) =\!\!\!\! \sum_{Y\subset M \atop \#Y={\rm
rank}(Y)=s}\!\!\!\!\! \alpha_Y E_{{c_{_Y}}}(x|Y),
 \end{equation}
 for all ${c}\in \C^n$ such that the denominators in $\alpha_Y$ do
not vanish.
\end{lemma}
Using Lemma \ref{le:E} we can now give an explicit formula for
$T(x|M)$.

\begin{theorem}\label{th:for}
 {\small
\begin{eqnarray}\label{eq:T}
 T(x|M) =
 \frac{1}{(n-s)!}\sum_{Y\subset M \atop \#Y={\rm rank}(Y)= s}\!\!\!\!
 {\alpha_Y}\, {|{\rm det}Y|^{-1}}
(-\theta_Y\cdot x)^{n-s} \chi_{{\rm
cone}(Y)}\left(x\right),\nonumber
 \end{eqnarray}}
  for all $c\in \C^n$ such that the denominators on the right-hand side do
not vanish, where both $\alpha_Y$  and $\theta_Y$ are defined in
Lemma \ref{le:E}.
\end{theorem}
\begin{proof} For an invertible $Y\subset M$, one has
$$
E_{c_Y}(x|Y) = \frac{1}{|\det Y|} \exp{(-\theta_Y\cdot x)}\chi_{{\rm
cone}(Y)}(x).
$$
 According to Lemma
\ref{le:E}, for $\rho\in \R\setminus \{0\}$, we have
\begin{eqnarray*}
& & E_{\rho {c}}(x|M) = \rho^{-n+s}\!\!\!\!\sum_{Y\subset M \atop
\#Y={\rm rank}(Y)= s}\!\!\! \alpha_Y  E_{\rho {c_{_Y}}}(x|Y).
 \end{eqnarray*}
Then the Taylor expansion of $E_{\rho c}(x|M)$  about $0$ in the
variable $\rho $ is
\begin{equation}\label{eq:taylor}
E_{\rho {c}} (x|M) = \sum_{l=0}^\infty
\rho^{l-n+s}p_l(x),
\end{equation}
where
$$
p_l(x) =\frac{1}{l!}\!\!\! \sum_{Y\subset M \atop
\#Y={\rm rank}(Y)= s}\!\! \alpha_Y |\det Y|^{-1} (-\theta_Y\cdot x)^l
\chi_{{\rm cone}(Y)}(x).
$$
The definition of $T(x|M)$ implies that it is the constant term in
(\ref{eq:taylor}). Hence, we have $T(x|M)=p_{n-s}(x)$. The theorem
follows.
\end{proof}

Brion's formula, which is obtained by Brion's Theorem, is useful for
computing the relative volumes of polytopes. We state it here.
\begin{theorem}{\rm (\cite{brion})}
Suppose that ${\mathcal P}$ is a simple rational convex
$d$-polytope. For a vertex cone  ${\mathcal K}_v$ of ${\mathcal P}$,
fix a set of generators $w_1(v),w_2(v),\ldots,w_d(v)\in \Z^d$. Then
$$
{\rm vol} ({\mathcal P})= \frac{\left(-1\right)^d}{d!}\sum_{{ v}\
\text{\rm a vertex of }
 \ {\mathcal P}}\frac{({ v}\cdot c)^d|\det(w_1(v),\ldots,w_d(v))|}{\prod_{k=1}^d(w_k(v)\cdot c)}
$$
for all $c\in \C^d$ such that the denominators on the right-hand
side do not vanish.
\end{theorem}
\begin{remark}
 Brion's formula requires the ${\mathcal V}$-representation and
the generators for each vertex cone, while the formula presented in
Theorem \ref{th:for} requires the ${\mathcal H}$-representation.
Hence, the formula in Theorem \ref{th:for} can be considered as a
dual version of Brion's.
\end{remark}

We next turn to another formula for computing the relative volume of
polytopes.
 We first introduce an iterative  formula for calculating
 the multivariate truncated power.

\begin{theorem}{\rm (\cite{numeircalT})}\label{Th:numT}
Let $M$ be an $s\times n$ matrix with columns $m_1, \ldots, m_n\in
\R^s\setminus \{ 0\}$ such that the origin is not contained in the interior of
${\rm conv}(M)$. Suppose that  $n>s+1$. For any $\lambda_1,\ldots,\lambda_n\in
\R$ and $x=\sum_{j=1}^n\lambda_j m_j$, we have
\begin{equation}\label{eq:recT}
T(x|M) = \frac{1}{n-s}\sum_{j=1}^n\lambda_j T(x|M
\setminus m_j).
\end{equation}
\end{theorem}

We next introduce Lasserre's  formula for  the volume of polytopes, which has
become well-known  today (see \cite{lasserre, volume1}). Consider the convex
polytope defined by
$$
D(b) := \{x\in \R^n |\,\, Ax\leq b \},
$$
where $A$ is an $s\times d$ matrix and $b$  an $s$-vector.
 The $i$th face of $D(b)$ is defined by
$$
D_i(b):= \{x\in \R^d|\,\, a_i\cdot x=b_i,\,\, Ax\leq b\},
$$
where $a_i$ is the $i$th row of $A$. We  set $V(d,A,b):={\rm
vol}_d(D(b))$ and set $V_i(d-1,A,b):={\rm vol}_d(D_i(b))$. Now we can
describe Lasserre's formula and show that it can be proved by
(\ref{eq:recT}).
\begin{theorem}\label{th:lasserre}{\rm(\cite{lasserre})}
If $V(d,A,b)$ is differentiable at $b$, then
\begin{equation}\label{eq:lasserre}
V(d,A,b)=\frac{1}{d}\sum_{i=1}^s \frac{b_i}{\|a_i\|}V_i(d-1,A,b).
\end{equation}

\end{theorem}
\begin{proof}
 Without loss of generality, we can suppose that
all points in $D(b)$ are non-negative, i.e., $D(b):=\{x\in
{\R}^d_+|\,\, Ax\leq b\}$. We first consider the case where each
entry in $A$ is an integer. By (\ref{Eq:truncatedvolume}) and
(\ref{Eq:truncated1}), when $A$ is an integer matrix,
$$
T(b|M)={\rm vol}(P)=V(d,A,b),
$$
where $ P:=\{x\in \R_+^{d+s}|\,\, Mx=b\}$ and $M:=(A,E_{s\times
s})$. Let  $e_i$ be the $s$-vector with $1$ at the $i$th position
and 0 for $j\neq i$. Using (\ref{eq:recT}), we obtain
\begin{equation}\label{eq:27}
 T(b|M) =\frac{1}{d}\,\sum_{i=1}^s b_i\, T(b|M\setminus e_i).
\end{equation}
Note that the $(d-1)$-polytope $D_i(b)$ lies in a hyperplane $\{x\in
\R^d|\,\, a_i\cdot x=b_i\}$ and that the $(n-1)$-dimensional volume
of the unit parallelogram in the hyperplane is $\frac{\|a_i\|}{{\rm
gcd}(a_{i1},\ldots,a_{id})}$ (cf. \cite{barv}), where $a_{ij}$ is
the $j$th entry in the vector $a_i$. So, one has
\begin{equation}\label{eq:volrat}
\frac{{\rm vol}_d(D_i(b))}{{\rm vol}(D_i(b))}=\frac{\|a_i\|}{{\rm
gcd}(a_{i1},\ldots,a_{id})}.
\end{equation}
By (\ref{Eq:truncated1}) and (\ref{eq:volrat}), we have
\begin{eqnarray*}
T(b|M\setminus e_i)&=&\frac{{\rm vol}(D_i(b))}{\# \{ [[M\setminus
e_i))\cap \Z^s\}}\,=\,\frac{{\rm vol}(D_i(b))}{{\rm
gcd}(a_{i1},\ldots,a_{id})}\\
&=&\frac{{\rm
vol}_d(D_i(b))}{\|a_i\|}\,=\,\frac{V_i(d-1,A,b)}{\|a_i\|}.
\end{eqnarray*}

Substituting $T(b|M)=V(d,A,b)$ and $T(b|M\setminus
e_i)=\frac{V_i(d-1,A,b)}{\|a_i\|}$ into (\ref{eq:27}), we get
(\ref{eq:lasserre}). By taking limit, (\ref{eq:lasserre}) holds for
any matrix $A$.
\end{proof}

\section{Integration of polynomials over polytopes}

In this section, we consider the problem of integrating of
polynomials over polytopes. Since each polynomial can be written as
the  sum of monomials, we only consider the monomial case.
  For every ${\bf k}=(k_1,\ldots,k_n)
\in \Z_+^{n}$ and $M=(m_1,\ldots,m_n)$ we set
$$
M^{\bf k}:=(
\underbrace{m_1,\ldots,m_1}_{k_1+1},\underbrace{m_2,\ldots,m_2}_{k_2+1},\ldots,\underbrace{m_n,\ldots,m_n}_{k_n+1}).
$$
The following theorem shows that the integration of  monomials can
be handled by the multivariate truncated power.
\begin{theorem}\label{th:lap} Suppose ${\bf k}=(k_1,\ldots,k_n)\in \Z_+^n$ and $f(u)=\prod_{j=1}^n
u_j^{k_j}$. Set $P:=\{u\in \R_+^n|\,\, Mu=x\}$. Then
$$
T(x|M^{\bf k}) = \frac{1}{{\bf k}!\cdot \sqrt{\,\det(MM^T)}}\int_{P}
f(u) du,
$$
where   ${\bf k}!:=k_1!\cdots k_n!$.
\end{theorem}
\begin{proof}
We set
$$
T_{\bf k}(x|M):= \frac{1}{\sqrt{\,\det(MM^T)}}\int_{P}f(u) du
$$
and consider its Laplace transform, i.e.,
$$
\widehat{T}_{\bf k}(\omega|M):= \int_{\R^s} \exp{(-\omega\cdot
x)}T_{\bf k}(x|M)\, dx .
$$
For every  $u$ in $P:=\{u\in \R_+^n|\,\, Mu=x\}$, we can write $u$
in a unique way as $u=z+y$ where $z\in {\rm ker}(M)$ and $y\bot {\rm
ker}(M)$. Note that $x=My$ and hence $\sqrt{\,\det(MM^T)}dy=dx$ (see
\cite{deboorbook}). Then we have
\begin{eqnarray*}
& &\int_{\R^s}\exp{(-\omega\cdot x)}T_{\bf k} (x |M)dx\\
&=&\frac{1}{\sqrt{\det(MM^T)}} \int_{\R^s}\exp{(-\omega\cdot x)}\int_{\{u\in \R_+^n|\,\, Mu=x\}} f(u)du dx\\
 &=&\int_{y\bot {\rm ker} M}\exp{(-\omega\cdot My)}\int_{z\in {\rm ker} M} \chi_{_{\R_+^n}}(z+y) f(z+y) dz
 dy\\
 &=&\int_{y \bot {\rm ker} M}\int_{z\in {\rm ker} M}\chi_{_{\R_+^n}}(z+y)  f(z+y)\exp{(-\omega\cdot M(z+y))} dz
 dy\\
  &=&\int _{{\R}_+^n}f(u)\exp{(-\omega\cdot Mu)}du = {\bf k}!\cdot\prod_{j=1}^n
\frac{1}{(\omega\cdot m_j)^{k_j+1}}.
\end{eqnarray*}
Also, note that
$$
 \widehat{T}(\omega|M^{\bf k}) = \prod_{j=1}^n \frac{1}{(\omega\cdot m_j)^{k_j+1}}.
$$
Hence the theorem follows from the inverse theorem for  Laplace
transform.

\end{proof}

\section{Multivariate box splines and the volume of cube slicing}

The {\it multivariate box spline} $B(\cdot |M)$ associated with $M$
is the distribution given by the rule (see \cite{deboor1,deboor2})
\begin{equation}\label{Eq:boxspline}
\int_{\R^s}\!\!\! B(x |M)\phi (x)dx = \int _{[0,1)^n}\!\! \phi
(Mu)du,\,\, \phi \in {\mathscr D}(\R^s).
\end{equation}
According to \cite{deboorbook}, one has
\begin{equation}\label{Eq:boxsplinevolume}
B(x|M) = \frac{{\rm vol}_n(P\cap
[0,1)^n)}{\sqrt{\,\det(MM^T)}},
\end{equation}
where $P:=\{y\in \R_+^n|\,\, My=x\}$. The formula
(\ref{Eq:boxsplinevolume}) shows the connection between the box
spline and the volume of cube slicing. Based on this connection, we
can study some interesting problems concerning the unit cube.

Recall that $Q_n$ is the unit cube in $\R^n$. We assume that  $H$ is an
$(n-1)$-dimensional hyperplane of $\R^n$ through the center of $Q_n$,
i.e.,
$$
H:=\{y\in \R^n|\,\, a_1y_1+\cdots+a_my_m=\sum_{i=1}^ma_i/2\},
$$
where  $1\leq m\leq n$ , $a_i$ is a  real number for $1\leq i\leq
m$.
 Based on the symmetry  of $Q_n$, we may assume that $a_i> 0$ for $1\leq i\leq m$. We set $A:=(a_1,\ldots,a_m)$.
Then from (\ref{Eq:boxsplinevolume}) we have
$$
{\rm vol}_n({H}\cap
Q_n)=\sqrt{\,\det(AA^T)}B\left({(a_1+\cdots+a_m)}/{2}|A\right).
$$
 By the symmetry  of the box spline, $B(x|A)$ achieves its
 maximum value at $({a_1+\cdots+a_m})/2$. So, as allowed to before,
Good's conjecture is equivalent to
\begin{equation}\label{eq:maxb}
 \max_x B(x|A) \geq  \frac{1}{\sqrt{\,\sum_{i=1}^m a_i^2}},
\end{equation}
 where  $1\leq m \leq n$. We next present a
 spline method for proving (\ref{eq:maxb}).

 \begin{theorem}
$$
\max_x B(x|(a_1,\ldots,a_m))\geq \frac{1}{\sqrt{\,\sum_{i=1}^m
a_i^2}},
$$
where $a_i$ are positive real numbers. The equality holds if and
only if $m=1$.
\end{theorem}
\begin{proof} Set
$$C(x|A):=B(x+{\sum_{i=1}^ma_i}/{2}|A).$$ The
Fourier transform of $C(x|A)$ is
$$
\widehat{C}(\omega|A)=\prod_{i=1}^m \frac{\sin(\omega a_i/2)}{\omega
a_i/2}.
$$
We can see that $\widehat{C}(0|A)=1$. According to the definition of
Fourier transform, we conclude that
\begin{equation}\label{eq:eq1}
\int_0^\infty \!\! t^2 C(t|A)dt = -\frac{1}{2} \widehat{C}''(0 |A)
=\frac{\sum_{i=1}^m a_i^2}{24}.
\end{equation}
Put $$S(t) := \int_0^t C(x|A)dx.$$
By
$$
\left(t\,\, {\max_x C(x|A)}\right)^2\,\, \geq\,\, S(t)^2
$$
we have
\begin{eqnarray}\label{eq:eq2}
\left({\max_x C(x|A)}\right)^2 \int_0^\infty\!\! t^2
C\left(t|A\right)dt &\geq& \int_0^\infty\!\!\!
S(t)^2 C\left(t|A\right)dt\\
&=&\frac{1}{3}\int_0^\infty\!\! \frac{dS(t)^3}{dt}dt\, =\,
\frac{1}{24},\nonumber
\end{eqnarray}
where the last equality follows from
$$\int_0^\infty
C(x|A)dx\,=\,\frac{\widehat{C}(0|A)}{2}\,=\,\frac{1}{2}.$$
 We combine (\ref{eq:eq1}) and (\ref{eq:eq2}) to obtain
$$
\max_xB(x|A)\,\, =\,\, \max_x C(x|A)\,\, \geq\,\,
\frac{1}{\sqrt{\,\sum_{i=1}^m a_i^2}}.
$$
From (\ref{eq:eq2}), we  see that the equality holds if and only
if $m=1$. \end{proof}

In the following and final theorem, we present an explicit formula for the
volume of $j$-slice of $Q_n$. The formula given in \cite{slice} is the
$(j=n-1)$-case in the next theorem.

\begin{theorem}
Suppose that $M$ is an $(n-j)\times n$ matrix with ${\rm
rank}(M)=n-j$ and let $P:=\{y\in \R_+^n |\,\, My=x\}$. Denote
$\Xi=\{0,1\}^n$ and $|\varepsilon|=\sum_{i=1}^n \varepsilon_i$. Then
we have
\begin{eqnarray}
& & {\rm vol}_n(P\cap Q_n)=\frac{\sqrt{\,\det(MM^T)}}{j!}\nonumber\\
& &\,\,\,\,\,\, \sum_{Y\subset M \atop \#Y={\rm rank}(Y)=s}\!\!\!
\alpha_Y |\det(Y)|^{-1} \sum_{\varepsilon\in
\Xi}(-1)^{|\varepsilon|}
 (-\theta_Y\cdot (x-M\varepsilon))^{j}  \chi_{{\rm
 cone}(Y)}(x-M\varepsilon),\nonumber
 \end{eqnarray}
 for all $c$ such that the denominators on the right-hand side
do not vanish, where both $\alpha_Y$ and $\theta_Y$ are defined in
Lemma \ref{le:E}.
\end{theorem}
\begin{proof}
By (\ref{Eq:boxsplinevolume}), one has
$$
{{\rm vol}_n(P\cap Q_n)}={\sqrt{\det(MM^T)}}B(x|M).
$$
Now we present an explicit formula for $B(x|M)$.
 Recall that $B(x|M)=\nabla_M T(x|M)$ (see \cite{deboorbook}), where $\nabla_{M}:=\prod_{i=1}^n\nabla_{m_i}$
 and $\nabla_{m_i}T(x|M):=T(x|M)-T(x-m_i|M)$.
From this formula, one obtains (see \cite{micchelli}):
\begin{equation}\label{eq:BT}
B(\cdot |M)=\sum_{\varepsilon\in \Xi}(-1)^{|\varepsilon |}T(\cdot
-M\varepsilon |M).
\end{equation}
The theorem is proved by putting (\ref{eq:T})
into (\ref{eq:BT}).
\end{proof}

\section{Examples}
\begin{example}
Set
$$D(z) := \{y\in \R_+^2|\,\, y_1+y_2\leq z;\,\, -2y_1+2y_2\leq
z;\,\, 2y_1-y_2\leq z\}.$$ The volume of $D(z)$ has been calculated
in \cite{lase}
 using Cauchy's Residue theorem. Here we can obtain it directly.
Based on (\ref{Eq:truncated1}), we have ${\rm vol}(D(z))=T({\bf
z}|M)$ where $$M=\begin{pmatrix}
  1 & 1 & 1 & 0 & 0 \\
  -2 & 2 & 0 & 1 & 0 \\
  2 & -1 & 0 & 0 & 1
\end{pmatrix},\,\,\,\, {\bf z}=(z,z,z)^T.
$$
We use $m_i$ to denote the $i$th column in $M$. A short calculation shows that
the cones spanned by the square matrices $(m_1,m_2,m_4)$, $(m_1,m_2,m_5)$,
$(m_1,m_3,m_4)$, $(m_2,m_3,m_5)$ and $(m_3,m_4,m_5)$ contain ${\bf z}$. We
select ${c}=(1,1,$ $1,1,1/2)$ in Theorem \ref{th:for} and obtain that $T({\bf
z}|M)=\frac{17}{48}z^2$ which agrees with the result presented in \cite{lase}.
\end{example}

\begin{example}
Set $\Omega_d:=\{{y}\in \R_+^d|\,\, \sum_{i=1}^d y_i \leq 1\}$. We
consider the problem of integrating of monomials over $\Omega_d$,
i.e.,
$$
J_d := \int_{\Omega_d}\!\! y_1^{k_1}\cdots y_d^{k_d}\,\, dy_1\cdots
dy_d.
$$
The value of $J_d$ is also calculated in \cite{simplex1,simplex2}.
Based on Theorem \ref{th:lap}, we can compute it easily. We set
${\rm e}_{d}:=(1,\ldots,1)\in \Z^d$. Using Theorem \ref{th:lap}, we
have
$$ J_d=k_1!\cdots k_d!\,\,\,
T(1|{\rm e}_{_{\sum_{i=1}^dk_i+d+1}}).
$$
It is well known that $T(x|{\rm e}_{d})=\frac{x_+^{d-1}}{(d-1)!}$.
Hence, we have
$$J_d=k_1!\cdots k_d!\,\, T(1|{\rm e}_{_{\sum_{i=1}^d
k_i+d+1}})=\frac{k_1!\cdots k_d!}{(\sum_{i=1}^d k_i+d)!}.$$
\end{example}

\vspace{1cm}

{\bf Acknowledgments.} The author is grateful to the referees for their
carefully reading of the manuscript and for their helpful comments on improving
the final version of this paper. The author also  thanks J. Gagelman,  O.~Holtz
and B.~Sturmfels  for  helpful discussions.

\end{document}